\documentclass[12pt]{article}
\usepackage{tikz}
\usetikzlibrary{intersections,calc,arrows}
\usepackage{amssymb}
\date{}
\makeatletter
\newcommand{\figcaption}[1]{\def\@captype{figure}\caption{#1}}
\newcommand{\tblcaption}[1]{\def\@captype{table}\caption{#1}}

\@addtoreset{equation}{section}
\makeatother
\newcommand{\qed}{\hbox{\rule[-2pt]{3pt}{6pt}}}

\begin{document}
\title {\bf Asymptotic behavior of bifurcation curves of nonlocal logistic equation of population dynamics}

\author{{\bf Tetsutaro Shibata}
\\
{\small Hiroshima University, Higashi-Hiroshima, 739-8527, Japan}
}

\maketitle
\footnote[0]{E-mail: tshibata@hiroshima-u.ac.jp
}
\footnote[0]{This work was supported by JSPS KAKENHI Grant Number JP25K07087.
}

\vspace{-0.5cm}

\begin{abstract}
We study the one-dimensional nonlocal Kirchhoff type bifurcation problem 
related to logistic equation of population dynamics.  
We establish the precise asymptotic formulas for bifurcation curve 
$\lambda = \lambda(\alpha)$ as $\alpha \to \infty$ in $L^2$-framework, 
where $\alpha:= \Vert u_\lambda \Vert_2$.

\end{abstract}


\vspace{0.5cm}


\section{Introduction} 		      

We consider the following one-dimensional nonlocal elliptic equation 
related to logistic equation of population dynamics

\begin{equation}
\left\{
\begin{array}{l}
-\left(a_1\Vert u\Vert_q^2+ a_2\Vert u\Vert_2^2 \right) u''(x)  + 
u(x)^p = \lambda u(x), \enskip 
x \in I:= (0,1),
\vspace{0.1cm}
\\
u(x) > 0, \enskip x\in I, 
\vspace{0.1cm}
\\
u(0) = u(1) = 0,
\end{array}
\right.
\end{equation}
where $a_1, a_2 \ge 0$ and $p, q > 1$ are given constants. 
We assume that $a_1 + a_2 > 0$. 
Further, 
$\lambda > 0$ is a bifurcation parameter.

Nonlocal elliptic problems have been investigated intensively by many authors 
and one of the main topics in this area is to study the existence,  nonexistence and the multiplicity of the solutions. We refer to [2-4, 6-13, 20] and the references therein. 
However, it seems that there are a few results which observe the nonlocal elliptic problems from a view point of bifurcation phenomena. 
In these studies, 
the bifurcation curves $\lambda$ were parameterized by 
$L^\infty$ norm of the solution $u_\lambda$ corresponding to 
$\lambda$ such as $\lambda = \lambda(\Vert u_\lambda\Vert_\infty)$ and the global structures 
of $\lambda(\Vert u_\lambda\Vert_\infty)$ have been investivated.  
We refer to [16-19, 21]. 

The purpose of this paper is to establish the precise asymptotic formulas for the bifurcation curves of the equation (1.1) in $L^2$-framework. 
That is, $\lambda$ is parameterized by 
$\alpha = \Vert u_\lambda\Vert_2$ such as $\lambda(\alpha)$ and 
consider the asymptotic behavior of $\lambda(\alpha)$ as 
$\alpha \to \infty$. 
As far as the author knows,  
there are few results 
to consider such nonlocal problem of logistic type as 
(1.1) from a view point of bifurcation problem in $L^2$-framework. 
In this sense, our results here are novel. 

To clarify our intension more precisely, 
we recall the following standard nonlinear eigenvalue problem of 
logistic type. Let 

\begin{equation}
\left\{
\begin{array}{l}
-w''(x) + w(x)^p = \gamma w(x), \enskip 
x \in I,
\vspace{0.1cm}
\\w(x) > 0, \enskip x\in I, 
\vspace{0.1cm}
\\
w(0) = w(1) = 0.
\end{array}
\right.
\end{equation}
Let $d > 0$ be an arbitrary given constant. Then we know from [1] that there exists a unique solution pair $(w_d, \gamma(d)) \in C^2(\bar{I}) \times \mathbb{R}_+$
of (1.2) with $\Vert w_d\Vert_2 = d$. Further, $\gamma$ is parameterized by $d$ such as $\gamma = \gamma(d)$, and 
it is called $L^2$-bifurcation curve. 
However, there are a few works to 
consider the precise global structure of $\gamma(d)$,
 since it is popular to investigate the global shape of the bifurcation curve 
$\gamma$ of (1.2) in $L^\infty$-
framework. Indeed, in many cases, $\gamma$ is parameterized by $L^\infty$ norm 
of the solution $w_\gamma$ associated with $\gamma$, namely, 
$\gamma = \gamma(\Vert w_\gamma \Vert_\infty)$.   
We emphasize 
that it is meaningful to treat the bifurcation problem (1.2) 
in $L^2$-framework, since the asymptotic behavior of 
$\gamma(\Vert w_\gamma \Vert_\infty)$ and $\gamma(d)$ are 
completely different from each other. It is well known (cf. [1]) that 
as $\Vert w_\gamma\Vert_\infty \to \infty$, then 
\begin{eqnarray}
\gamma(\Vert w_\gamma \Vert_\infty) = \Vert w_\gamma\Vert_\infty^{p-1} + O(1).
\end{eqnarray} 
On the other hand, it was shown in [14] that, for $d \gg 1$, 
\begin{eqnarray}
\gamma(d) &=& d^{p-1} + C_1d^{(p-1)/2} + O(1),
\end{eqnarray}
where
\begin{eqnarray}
C_1 = (p+3)\int_0^1 \sqrt{\frac{p-1}{p+1} - s^2 + \frac{2}{p+1}s^{p+1}}ds.
\end{eqnarray}
It is clear that (1.3) is affected only the behavior of $w_\gamma$ at the center of the interval $I$. On the other hand, (1.4) is affected 
not only the shape of $w_d$ in the interior of $I$ but also the behavior of $w_d$ 
near the boundary of $I$. Indeed, the second term of (1.4) comes from the asymptotic 
behavior of the slope of $w_d$ near the boundary of $I$. 

Motivated by the result mentioned above, we here concentrate on the 
effect of the nonlocal terms to the asymptotics of the 
equation and establish the asymptotic formulas for $\lambda(\alpha)$ as $\alpha \to \infty$, 
which are different from (1.4). By these results, we understand well how the nonlocal terms give effect to the 
asymptotic behavior of $\lambda(\alpha)$ as $\alpha \to \infty$. 

Now we state our results.

\vspace{0.2cm}

\noindent
{\bf Theorem 1.1.} {\it Assume that $p > 3$. Let 
$a_1, a_2 \ge 0$ be constants satisfying   
$a_1 + a_2 > 0$. Then for any given constant $\alpha > 0$, there exists a unique solution pair 
 $(u_\alpha, \lambda(\alpha)) \in C^2(\bar{I}) \times \mathbb{R}_+$ of (1.1) satisfying $\Vert u_\alpha \Vert = \alpha$. Furthermore, as $\alpha \to \infty$, 
\begin{eqnarray}
\lambda(\alpha) &=&\alpha^{p-1}\left\{1 + C_1(a_1+ a_2)^{1/2} \alpha^{-(p-3)/2}
+ O(\alpha^{-(p-3)})\right\}.
\end{eqnarray}
}

\vspace{0.2cm}

\noindent
{\bf Theorem 1.2.} {\it Assume that $p = 3$. Let $a_1, a_2 \ge 0$ 
be constants satisfying  $a_1+ a_2 > 0$. 
Then for any given constant $\alpha > 0$, there exists a unique solution pair $(u_\alpha, \lambda(\alpha)) \in C^2(\bar{I}) \times \mathbb{R}_+$ satisfies (1.1), which is represented as follows. 
There exists a unique constant $d_1 > 0$, $\gamma(d_1) > 0$ and 
$h_{d_1} = \alpha/d_{1}$ such that 
$(w_{d_1}, \gamma(d_1))$ satisfies (1.2) with $\Vert w_{d_1}\Vert_2 = d_1$ and 
$(u_\alpha, \lambda(\alpha)) = (h_{d_1}w_{d_1}, h_{d_1}^2\gamma(d_1))$.
}

\vspace{0.2cm}

\noindent
{\bf Theorem 1.3.} {\it Assume that $1 < p < 3$. Let $a_1, a_2 \ge 0$ 
be constants 
satisfying $a_1+ a_2 > 0$. 
Let $\alpha > 0$ be a given constant. Then as $\alpha \to \infty$
\begin{eqnarray}
\lambda(\alpha) &=& \pi^{2\{q(p-3)-(p-1)\}/(q(p-3))}E_1E_3^{-1}\alpha^2 
\\
&& \qquad \qquad \times 
\left[
1 + \left\{\frac{E_2}{E_1} + E_4 - E_3^{(p-3)/2}E_5\right\}E_3^{-(p-3)/2} \alpha^{p-3} 
+ o(\alpha^{p-3})
\right],
\nonumber
\end{eqnarray}
where 
\begin{eqnarray}
E_1 &:=& a_12^{(q+2)/q}A_1^{2/q} + a_2\pi^{2/q},
\\
E_2 &:=& a_12^{(q+2)/q}A_1^{2/q}\left(A_4+ \frac{2}{q}\frac{A_2}{A_1}\right) 
+ \frac{2a_2\pi^{(2-q)/q}}{q}A_3,
\\
E_3 &:=& \pi^{-4/((p-1)q)}E_1^{2/(p-3)}, 
\\
E_4 &:=& 2\frac{q(p-3)-(p-1)}{q(p-3)\pi}A_3,
\\
E_5 &:=& \left(\frac{2}{p-3}\frac{E_2}{E_1} - A_6\right)
\pi^{2(p-3)/((p-1)q)}E_1^{-1},
\\
A_1 &:=& \int_0^1 \frac{s^q}{\sqrt{1-s^2}}ds,
\\
A_2 &:=& \frac{(\sqrt{2})^{p-1}}{(p+1)\pi^2}\int_0^1 \frac{s^q(1-s^{p+1})}{(1-s^2)^{3/2}}
ds,
\\
A_3 &:=& \frac{2}{(p+1)\pi^2}(\sqrt{2})^{p-1}\int_0^1 \frac{1-s^{p+1}}{(1-s^2)^{3/2}}ds,
\\
A_4 &:=& \frac{1}{\pi}\left(A_3 - 4A_2\right),
\\
A_5&:=& 
\frac{1}{(p+1)\pi^2}\int_0^1 \frac{s^2(1-s^{p+1})}{(1-s^2)^{3/2}}ds,
\\
A_6&:=& \frac{4}{(p-3)q\pi}A_3.
\end{eqnarray}
}

\vspace{0.2cm}

The remainder of this paper is organized as follows. In Section 2, 
we prove Theorem 1.1 with the aid of the results in [14, 15]. 
In Section 3, we prove Theorems 1.2 and 1.3 by Taylor 
expansion and complicated direct calculation.

\section{Proof of Theorem 1.1} 		      

In what follows, we use the notations defined in Section 1. 
We begin with the existence of the solution pair 
$(u_\alpha, \lambda(\alpha))$.

\vspace{0.2cm}

\noindent
{\bf Lemma 2.1.} 
{\it Let $\alpha > 0$ be a fixed constant. Then there exists a 
unique solution pair 
$(u_\alpha, \lambda(\alpha))$ of (1.1) with $u_\alpha = hw_d$ for 
some $h > 0$ and $d > 0$. 
}

\vspace{0.2cm}

\noindent
{\it Proof.} 
Assume that $u_\alpha$ is a solution of (1.1) with $\lambda = \lambda(\alpha)$. We put
\begin{eqnarray}
\beta:= \beta(\alpha) = a_1\Vert u_\alpha\Vert_q^2 + a_2\Vert u_\alpha\Vert_2^2.
\end{eqnarray}
Then we have 
\begin{eqnarray}
-\beta u_\alpha'' + u_\alpha^p = \lambda u_\alpha.
\end{eqnarray}
For constants $h, d > 0$, we put $w_d:= h^{-1}u_\alpha$. 
Since $u_\alpha = hw_d$, by (2.2), we have 
\begin{eqnarray}
-\beta h w_d'' + h^p w_d^p= \lambda(\alpha) hw_d.
\end{eqnarray}
Namely, 
\begin{eqnarray}
-w_d'' + \frac{h^{p-1}}{\beta} w_d^p = \frac{\lambda(\alpha)}{\beta} w_d.
\end{eqnarray}
Let $h$ satisfy 
\begin{eqnarray}
h^{p-1} = \beta =  a_1h^2\Vert w_d\Vert_q^2 + a_2h^2\Vert w_d\Vert_2^2,
\end{eqnarray}
namely
\begin{eqnarray}
h =  (a_1\Vert w_d\Vert_q^2 + a_2\Vert w_d\Vert_2^2)^{1/(p-3)},
\end{eqnarray}
then we see from [14] that $(w_d, \gamma(d)) 
= (w_d, \frac{\lambda(\alpha)}{\beta})$ satisfies (1.2). Moreover, 
\begin{eqnarray}
\alpha = \Vert u_\alpha\Vert_2 = h\Vert w_d\Vert_2 =  (a_1\Vert w_d\Vert_q^2 + a_2\Vert w_d\Vert_2^2)^{1/(p-3)}\Vert w_d\Vert_2:= g(d).
\end{eqnarray}
We know from [14] that if $0 < d_1 < d_2$, then $w_{d_1} < w_{d_2}$ for $0 < x < 1$. Therefore, we see that $g(d)$ is strictly increasing function of $d$ and $g(d) \to 0$ as $d \to 0$. 
This implies that $d$ is a strictly increasing function of $\alpha > 0$, 
namely, $d = d_\alpha = g^{-1}(\alpha)$. 
Namely, there exists a unique $d_\alpha > 0$ such that $\alpha = g(d_\alpha)$ 
for any given $\alpha > 0$. By (2.1), we know that $\beta$ is a 
function of $\alpha$. Then $\lambda$ is determined uniquely 
by $\alpha$ such as 
$\lambda(\alpha) = \beta(\alpha)\gamma(g^{-1}(\alpha)) 
= \beta(\alpha))
\gamma(d_\alpha)$. 
We understand from (2.7) that $(u_\alpha, \lambda(\alpha))$ 
and $(w_{d_\alpha}, \gamma(d_\alpha))$ is one to one correspondence. 
This implies the unique existence of $(u_\alpha, \lambda(\alpha))$ 
for a given $\alpha$. Thus the proof is complete. \qed

\vspace{0.2cm}

\noindent
{\bf Proof of Theorem 1.1.} 
We know from [15, Proposition 2.1] that for $d \gg 1$,  
\begin{eqnarray}
\Vert w_d\Vert_q^{p-1} &=& \gamma(d)\left(1 - \frac{C(q)}{\sqrt{\gamma(d)}}\right)
^{(p-1)/q} + O(\gamma(d)e^{-\delta_1\gamma(d)}),
\end{eqnarray}
where $\delta_1 > 0$ is a constant. Here, 
\begin{eqnarray}
{C(q)} := 2\int_0^1 \frac{1-s^q}{\sqrt{1-s^2- \frac{2}{p+1}(1-s^{p+1})}}ds.
\end{eqnarray}
By (1.4), (2.8) and Taylor expansion, we have 
\begin{eqnarray}
\Vert w_d\Vert_q^2 &=& \gamma(d)^{2/(p-1)}
\left(1-\frac{2}{q}\frac{C(d)}{\sqrt{\gamma(d)}} + O(\gamma^{-1})\right)
\\
&=& \left(d^{p-1} + C_1d^{(p-1)/2} + O(1)\right)^{2/(p-1)}
\left(1-\frac{2}{q}\frac{C(d)}{\sqrt{\gamma(d)}} + O(\gamma^{-1})\right)
\nonumber
\\
&=& d^2\left(1 + \frac{2}{p-1}C_1d^{-(p-1)/2} + O(d^{-(p-1)})\right)
\nonumber
\\
&&\times
\left\{1 - \frac{2}{q}C(q)
d^{-(p-1)/2} + O(d^{-(p-1)})
\right\}
\nonumber
\\
&=& d^2\left\{1 + \left(\frac{2}{p-1}C_1-\frac{2}{q}C(q)\right)d^{-(p-1)/2} 
+ O(d^{-(p-1)})\right\}
\nonumber
\\
&=:& d^2D(d).
\nonumber
\end{eqnarray}
By this and (2.7), 
we have 
\begin{eqnarray}
d = \alpha^{(p-3)/(p-1)}(a_1D(d) + a_2)^{-1/(p-1)}.
\end{eqnarray}
By this and (2.10), we have 
\begin{eqnarray}
D(d)&=& 1 + \left(\frac{2}{p-1}C_1 - \frac{2}{q}C(q)\right)d^{-(p-1)/2} + O(d^{-(p-1)})
\\
&=& 1 + \left(\frac{2}{p-1}C_1 - \frac{2}{q}C(q)\right)
\left\{\alpha^{(p-3)/(p-1)}(a_1D(d) + a_2)^{-(1/(p-1)}\right\}^{-(p-1)/2}
\nonumber
\\
&& + O(\alpha^{-(p-3)})
\nonumber
\\
&=& 1 + \left(\frac{2}{p-1}C_1 - \frac{2}{q}C(q)\right)\alpha^{-(p-3)/2}(a_1D(d)+a_2)^{1/2}
+ O(\alpha^{-(p-3)})
\nonumber
\\
 &=& 1 + \left(\frac{2}{p-1}C_1 - \frac{2}{q}C(q)\right)\alpha^{-(p-3)/2}(a_1+a_2)^{1/2}
+ O(\alpha^{-(p-3)}).
\nonumber
\end{eqnarray}
This along with (1.4), (2.6) and (2.10) implies that
\begin{eqnarray}
\lambda(\alpha) &=& \beta(d)\gamma(d) = h(d)^{p-1}\gamma(d) 
\\
&=& 
(a_1\Vert w_d\Vert_q^2 
+ a_2d^2)^{(p-1)/(p-3)}\gamma(d)
\nonumber
\\
&=& (a_1d^2D(d) + a_2d^2)^{(p-1)/(p-3)}d^{p-1}(1 + C_1d^{-(p-1)/2} + O(d^{-(p-1)})
\nonumber
\\
&=& d^{(p-1)^2/(p-3)}(a_1D(d) + a_2)^{(p-1)/(p-3)}(1 + C_1d^{-(p-1)/2} + O(d^{-(p-1)}).
\nonumber
\end{eqnarray}
By this and (2.11), we have 
\begin{eqnarray}
\lambda(\alpha) &=& \alpha^{p-1}\left\{1 + C_1d^{-(p-1)/2} + O(d^{-(p-1)})\right\}
\\
&=& \alpha^{p-1}\left\{1 + C_1\alpha^{-(p-3)/2}(a_1D(d) + a_2)^{1/2} + O(\alpha^{-(p-3)})
\right\}
\nonumber
\\
&=&\alpha^{p-1}\left\{1 + C_1\alpha^{-(p-3)/2}(a_1+ a_2)^{1/2} + O(\alpha^{-(p-3)})
\right\}.
\nonumber
\end{eqnarray}
This implies (1.3). Thus the proof is complete. \qed

\section{Proof of Theorems 1.2 and 1.3}

We begin with the proof of Theorem 1.2, which is the same 
argument 
as that used in the proof of Lemma 2.1.

\vspace{0.2cm}

\noindent
{\bf Proof of Theorem 1.2.} Let $p = 3$. We apply the argument in the previous section to the case $p = 3$. Then by (2.5), we have 
\begin{eqnarray}
1 = a_1\Vert w_d\Vert_q^2 + a_2\Vert w_d\Vert_2^2.
\end{eqnarray}
Since the r.h.s. of (3.1) is strictly increasing function of $d > 0$ and 
tends to $0$ as $d \to 0$, there exists a unique constant $d = d_1> 0$ and $(w_{d_1}, \gamma(d_1)) \in C^2(\bar I) \times \mathbb{R}_+$ satisfying (1.2) and (3.1) with 
$\Vert w_{d_1}\Vert_2 = d_1$. By (2.5), we have $\beta = h^2$ and 
$\frac{\lambda}{\beta} = \gamma(d_1)$. Further, $\alpha = \Vert u_\alpha\Vert_2 
= h\Vert w_{d_1}\Vert_2 = hd_1$. By this, we have 
\begin{eqnarray}
\lambda(\alpha) &=& \beta\gamma(d_1) = h^2\gamma(d_1) = \frac{\alpha^2}{d_1^2}
\gamma(d_1).
\end{eqnarray}
Thus the proof of Theorem 1.2 is complete. \qed

\vspace{0.2cm}

Now we prove Theorem 1.3. For 
an arbitrary given constant $\alpha > 0$, the proof of the unique 
existence of the solution pair $(u_\alpha, \lambda(\alpha))$ 
of (1.1) with $\Vert u_\alpha \Vert_2 = \alpha$ is the same as 
that of Lemma 2.1. We also find that 
Lemma 2.1 is also true for the case $1 < p < 3$. So 
we use the same notations as those defined in the proof of Lemma 2.1 
in what follows. 

\vspace{0.2cm}

\noindent
{\bf Lemma 3.1.} {\it Let $1 < p < 3$. Then $d \to 0$ as $\alpha 
\to \infty$.
}

\vspace{0.2cm}

\noindent
{\it Proof.} 
We put $u = hw_d$ and $\Vert w_d\Vert_2 = d$. Then we have 
$\alpha = hd$. We first 
assume that there exists a constant $M > 0$ such that 
$M^{-1} < d < M$. 
Then we see from [1] that $\Vert w_d\Vert_q$ is bounded.  
Then by (2.5), we see that $h$ is bounded. Then $\alpha = hd$ 
is bounded. This is a contradiction. 
Next, we assume that $d \to \infty$ as $\alpha \to \infty$. 
Then by [1, 14], we know that $w_d(x) = d(1 + o(1))$ for $x \in I$. 
By this, (1.3), (1.4) and (2.6), we see that $h \sim d^{2/(p-3)}$. By this and (2.7), we have 
$\alpha \sim d^{(p-1)/(p-3)} \to 0$. This is a contradiction. 
Therefore, $d \to 0$ as $\alpha \to \infty$. Thus the proof is complete. \qed

\vspace{0.2cm}

By Lemma 3.1, let $0 < d \ll 1$ in what follows.  
By Lemma 3.1 and (1.2), we see that $w_d(x) \to \sqrt{2}d \sin \pi x$ in $C^1(\bar{I})$ as 
$d \to 0$, since $p > 1$. Moreover, 
$\gamma(d) \to \gamma(0) = \pi^2$ as $d \to 0$, where 
$\pi^2$ is the first eigenvalue of the 
linear eigenvalue problem corresponding to (1.2). Recall that 
by (2.13), we know   
\begin{eqnarray}
\lambda(\alpha) &=& \beta\gamma(d) = h^{p-1}\gamma(d)
\\
&=& (a_1\Vert w_d\Vert_q^2 + a_2\Vert w_d\Vert_2^2)^{(p-1)/(p-3)}\gamma(d).
\nonumber
\end{eqnarray} 
We calculate $\Vert w_d\Vert_q$ by using the time map method in [16]. 
For simplicity, we write $w = w_d$, $k:=\Vert w_d\Vert_\infty = \sqrt{2}d(1 + o(1))$ and $\gamma = \gamma(d) = \pi^2(1 + o(1))$ 
in what follows. 
\vspace{0.2cm}

\noindent
{\bf Lemma 3.2.} {\it As $d \to 0$,
\begin{eqnarray}
\Vert w\Vert_q^2 &=& \frac{(2A_1)^{2/q}k^2}{\gamma^{1/q}}\left( 1 + \frac{2}{q}\frac{A_2}{A_1}d^{p-1} + o(d^{p-1})\right).
\end{eqnarray}
}
{\it Proof.}  
If $w$ satisfies (1.2), then by [5], we have 
\begin{eqnarray}
w(x) &=& w(1-x), \quad 0 \le x \le 1,
\\
w'(x) &>& 0, \quad 0 < x < \frac12,
\\
\Vert w\Vert_\infty &:=& \max_{x \in I} w(x) = w\left(\frac12\right).
\end{eqnarray}
By (1.2), for 
$x \in \bar{I}$, we have 
\begin{eqnarray}
(w''(x) + \gamma w(x) - w(x)^p)w(x) = 0.
\end{eqnarray}
By this, (3.7) and putting $x = 1/2$, we have 
\begin{eqnarray}
\frac12(w'(x))^2 + \frac12 \gamma w(x)^2 - \frac{1}{p+1}w(x)^{p+1} &=& \mbox{constant} 
\\
&=& \frac12\gamma k^2 - \frac{1}{p+1}k^{p+1}.
\nonumber
\end{eqnarray}
By this and (3.6), for $0 \le x \le 1/2$, we have
\begin{eqnarray}
w'(x) &=& \sqrt{\gamma (k^2 - w(x)^2) - \frac{2}{p+1}(k^{p+1} - w(x)^{p+1})}
\end{eqnarray} 
By this, (3.5) and putting $\theta = ks = w(x)$, we have 
\begin{eqnarray}
\Vert w\Vert_q^q &=& 2\int_0^{1/2} \frac{w(x)^qw'(x)}
{\sqrt{\gamma (k^2 - w(x)^2) - \frac{2}{p+1}(k^{p+1} - w(x)^{p+1})}}dx
\\
&=& 2\int_0^k \frac{\theta^q}{\sqrt{\gamma (k^2 - \theta^2) 
- \frac{2}{p+1}(k^{p+1} - \theta^{p+1})}}d\theta
\nonumber
\\
&=& 2k^q\int_0^1 \frac{s^q}{\sqrt{\gamma(1-s^2) - \frac{2}{p+1}k^{p-1}(1-s^{p+1})}}ds
\nonumber
\\
&=& \frac{2k^q}{\sqrt{\gamma}}\int_0^1 \frac{s^q}{\sqrt{(1-s^2) - \frac{2}{p+1}
\frac{k^{p-1}}{\gamma}(1-s^{p+1})}}ds.
\nonumber
\end{eqnarray}
Since $k = \sqrt{2}d(1 + o(1))$ for $0 < d \ll 1$, by (3.11) and Taylor expansion, we have 
\begin{eqnarray}
\Vert w\Vert_q^q &=& \frac{2k^q}{\sqrt{\gamma}}\int_0^1 \frac{s^q}{\sqrt{1-s^2}}
\left\{1 + \frac{1}{(p+1)\pi^2}k^{p-1}
\frac{1-s^{p+1}}{1-s^2}(1 + o(1))\right\}ds
\\
&=& \frac{2k^q}{\sqrt{\gamma}}\left\{A_1 + 
A_2d^{p-1} + o(d^{p-1})\right\}.
\nonumber
\end{eqnarray}
By this and Taylor expansion, we have  
\begin{eqnarray}
\Vert w\Vert_q^2 &=& \frac{(2A_1)^{2/q}k^2}{\gamma^{1/q}}\left( 1 + \frac{2}{q}\frac{A_2}{A_1}d^{p-1} + o(d^{p-1})\right).
\end{eqnarray}
This implies our conclusion. Thus the proof is complete. \qed

\vspace{0.2cm}

We next calculate $\gamma$ precisely. 

\vspace{0.2cm}

\noindent
{\bf Lemma 3.3.} {\it As $d \to 0$,
\begin{eqnarray}
\sqrt{\gamma} &=& \pi + A_3d^{p-1} + o(d^{p-1}),
\\
\gamma &=& \pi^2 + 2\pi A_3d^{p-1} + o(d^{p-1}).
\end{eqnarray}
}
{\it Proof.} By (3.10) and Taylor expansion, we have 
\begin{eqnarray}
\frac12 &=& \frac12 \int_0^1 dx = \int_0^{1/2}\frac{w'(x)}
{\sqrt{\gamma (k^2 - w(x)^2) - \frac{2}{p+1}(k^{p+1} - w(x)^{p+1})}}dx
\\
&=& 
 \frac{1}{\sqrt{\gamma}}\int_0^1 \frac{1}{\sqrt{1-s^2}}
\left\{1 + \frac{1}{(p+1)\gamma}k^{p-1}\frac{1-s^{p+1}}{1-s^2}(1 + o(1))\right\}ds
\nonumber
\\
&=& \frac{1}{\sqrt{\gamma}}\left\{ \frac{\pi}{2} + \frac{1}{(p+1)\pi^2}
(\sqrt{2}d)^{p-1}\int_0^1 \frac{1-s^{p+1}}{(1-s^2)^{3/2}}(1 + o(1))\right\}ds.
\nonumber
\end{eqnarray}
By this, we have 
\begin{eqnarray}
\sqrt{\gamma} &=& \pi + A_3d^{p-1} + o(d^{p-1}).
\end{eqnarray}
This implies (3.14). (3.15) follows immediately from (3.14). \qed

\vspace{0.2cm}

We now calculate $k^2$ precisely. 

\vspace{0.2cm}

\noindent
{\bf Lemma 3.4.} {\it As $d \to 0$, 
\begin{eqnarray}
k^2 &=& 2d^{2}\left\{1 + A_4d^{p-1} 
+ o(d^{p-1})\right\}.
\end{eqnarray}
}
{\it Proof.} We note that $A_1 = \pi/4$ when $q = 2$. By putting $q = 2$ in (3.12), we have 
\begin{eqnarray}
d^2 &=& \Vert w\Vert_2^2 = \frac{2k^2}{\sqrt{\gamma}}
\left\{\frac{\pi}{4} + A_2d^{p-1} + o(d^{p-1})\right\}.
\end{eqnarray}
This implies that 
\begin{eqnarray}
\frac{\sqrt{\gamma}d^2}{2} &=& k^2\left\{\frac{\pi}{4} +
A_2 d^{p-1} + o(d^{p-1})\right\}.
\end{eqnarray} 
By this, (3.14) and Taylor expansion, we have 
\begin{eqnarray}
k^2 &=& \frac{d^2}{2}\frac{\left\{ \pi + A_3d^{p-1} + o(d^{p-1})\right\}}
{\left\{\frac{\pi}{4} +
A_2 d^{p-1} + o(d^{p-1})\right\}}
\\
&=& 2d^2\left\{1 + \frac{1}{\pi}A_3d^{p-1} + o(d^{p-1})\right\}
\left\{1 - \frac{4}{\pi}A_2 d^{p-1} + o(d^{p-1})\right\}
\nonumber
\\
&=& 2d^{2}\left\{1 + A_4d^{p-1} 
+ o(d^{p-1})\right\}.
\nonumber
\end{eqnarray}
This implies (3.18). Thus the proof is complete. \qed

\vspace{0.2cm}

Now we represent $d$ by using $\alpha$ precisely. 

\vspace{0.2cm}

\noindent
{\bf Lemma 3.5.} {\it As $d \to 0$
\begin{eqnarray}
d^{p-1} = E_3^{-(p-3)/2}\alpha^{p-3}
\left(1 - \frac{p-3}{2}E_5\alpha^{p-3} + o(\alpha^{p-3})\right).
\end{eqnarray}
}
{\it Proof.}  
By (2.6), Lemmas 3.2 and 3.4, we have
\begin{eqnarray}
\alpha^2 &=& h^2d^2 = \left(a_1\frac{(2A_1)^{2/q}k^2}{\gamma^{1/q}}\left( 1 + \frac{2}{q}\frac{A_2}{A_1}d^{p-1} + o(d^{p-1})\right) + a_2d^2
\right)^{2/(p-3)} d^{2}
\\
&=& \left\{a_1\frac{(2A_1)^{2/q}2(1 + A_4d^{p-1} + o(d^{p-1}))}{\gamma^{1/q}}\left( 1 + \frac{2}{q}\frac{A_2}{A_1}d^{p-1} + o(d^{p-1})\right) + a_2\right\}^{(2/(p-3)}
\nonumber
\\
&&\times d^{2(p-1)/(p-3)}
\nonumber
\\
&=& \left\{a_12^{(q+2)/q}A_1^{2/q}(1 + A_4d^{p-1} + o(d^{p-1}))
(1 + \frac{2}{q}\frac{A_2}{A_1}d^{p-1} + o(d^{p-1})) 
+ a_2\gamma^{1/q}\right\}^{2/(p-3)}
\nonumber
\\
&&\times \left\{\pi^2 + 2\pi A_3d^{p-1} 
+ o(d^{p-1})\right\}^{-2/((p-3)q)}d^{2(p-1)/(p-3)}
\nonumber
\\
&=& \left\{a_12^{(q+2)/q}A_1^{2/q}\left(1 + \left(A_4+\frac{2}{q}
\frac{A_2}{A_1}\right)d^{p-1} + o(d^{p-1})\right) 
\right.
\nonumber
\\
&&\left. \qquad \qquad \qquad \qquad \qquad 
+ a_2\left(\pi^2 + 2\pi A_3d^{p-1} + o(d^{p-1})\right)^{1/q}
\right\}^{2/(p-3)}
\nonumber
\\
&&\times \pi^{-4/((p-3)q)}\left\{
1 - \frac{4}{(p-3)q\pi}A_3d^{p-1} + o(d^{p-1})\right\}
d^{2(p-1)/(p-3)}
\nonumber
\\
&=& \pi^{-4/((p-3)q)}
\left[(a_12^{(q+2)/q}A_1^{2/q} + a_2\pi^{2/q}) 
\right.
\nonumber
\\
&&
\left.+
\left\{a_12^{(q+2)/q}A_1^{2/q}\left(A_4 + \frac{2}{q}\frac{A_2}{A_1}\right) 
+ \frac{2a_2\pi^{(2-q)/q}}{q}A_3\right\}d^{p-1} + o(d^{p-1}) 
\right]^{2/(p-3)}
\nonumber
\\
&&
\qquad \qquad \quad 
\times
\left\{1 - \frac{4}{(p-3)q\pi}A_3d^{p-1} 
+ o(d^{p-1})\right\}d^{2(p-1)/(p-3)}
\nonumber
\\
&=&\pi^{-4/((p-1)q)}(E_1 + E_2d^{p-1} + o(d^{p-1}))^{2/(p-3)}
\nonumber
\\
&& \qquad \qquad \quad \qquad \qquad \times
(1 - A_6d^{p-1} + o(d^{p-1}))d^{2(p-1)/(p-3)}.
\nonumber 
\end{eqnarray} 
By this, we have $d^{p-1} = \pi^{2(p-3)/((p-1)q)}E_1^{-1}\alpha^{p-3}(1 + o(1))$. 
By this, (3.23) and Taylor expansion, we have 
\begin{eqnarray}
\alpha^2 &=& E_3\left\{1 + \frac{2}{p-3}\frac{E_2}{E_1}d^{p-1} + o(d^{p-1})\right\}
\left\{1 - A_6d^{p-1} + o(d^{p-1})\right\}
\\
&& \qquad \times d^{2(p-1)/(p-3)}
\nonumber
\\
&=& E_3\left\{1 + \left(\frac{2}{p-3}\frac{E_2}{E_1} - A_6\right)
d^{p-1} + o(d^{p-1})\right\}d^{2(p-1)/(p-3)}
\nonumber
\\
&=& E_3\left\{1 + \left(\frac{2}{p-3}\frac{E_2}{E_1} - A_6\right)
\pi^{2(p-3)/((p-1)q)}E_1^{-1}\alpha^{p-3} + o(\alpha^{p-3})\right\}
\nonumber
\\
&& \qquad \times 
d^{2(p-1)/(p-3)}
\nonumber
\\
&=& E_3\left\{1 + E_5\alpha^{p-3} + o(\alpha^{p-3})\right\}
d^{2(p-1)/(p-3)}.
\nonumber
\end{eqnarray}
By this, we have 
\begin{eqnarray}
d^{p-1} = E_3^{-(p-3)/2}\alpha^{p-3}
\left(1 - \frac{p-3}{2}E_5\alpha^{p-3} + o(\alpha^{p-3})\right).
\end{eqnarray}
Thus the proof is complete. \qed

\vspace{0.2cm}

\noindent
{\bf Proof of Theorem 1.3.} By (3.3), Lemmas 3.2, 3.4 and 3.5, we have 
\begin{eqnarray}
\lambda(\alpha) &=& h^{p-1}\gamma 
= (a_1\Vert w\Vert_q^2 + a_2d^2)^{(p-1)/(p-3)}\gamma
\\
&=& \left\{a_1\frac{(2A_1)^{2/q}k^2}{\gamma^{1/q}}\left( 1 + \frac{2}{q}\frac{A_2}{A_1}d^{p-1} + o(d^{p-1})\right) + a_2d^2\right\}^{(p-1)/(p-3)} \gamma
\nonumber
\\
&=&\left\{a_1(2A_1)^{2/q}2d^{2}(1 + A_4d^{p-1} + o(d^{p-1}))\left( 1 + \frac{2}{q}\frac{A_2}{A_1}d^{p-1} + o(d^{p-1})\right) 
\right.
\nonumber
\\
&& \left. \qquad 
+ a_2d^2\gamma^{1/q}
\right\}^{(p-1)/(p-3)} \left\{
\pi^2 + 2\pi A_3d^{p-1} + o(d^{p-1})
\right\}^{\{q(p-3)-(p-1)\}/(q(p-3))}
\nonumber
\\
&=& \left\{a_12^{(q+2)/q}A_1^{2/q}(1 + A_4d^{p-1} + o(d^{p-1}))\left( 1 + \frac{2}{q}\frac{A_2}{A_1}d^{p-1} + o(d^{p-1})\right) 
\right.
\nonumber
\\
&& \left. \qquad 
+ a_2\gamma^{1/q}\right\}^{(p-1)/(p-3)} 
\nonumber
\\
&&\times d^{(2(p-1)/(p-3)}\left\{
\pi^2 + 2\pi A_3d^{p-1} + o(d^{p-1})
\right\}^{{\{q(p-3)-(p-1)\}/(q(p-3))}}.
\nonumber
\end{eqnarray}
By this, we have 
\begin{eqnarray}
\lambda(\alpha) &=& \left\{a_12^{(q+2)/q}A_1^{2/q}(1 + A_4d^{p-1} + o(d^{p-1}))
\left( 1 + \frac{2}{q}\frac{A_2}{A_1}d^{p-1} + o(d^{p-1})\right) 
\right.
\\
&& \left. \qquad \qquad + a_2\left(\pi^2 + 2\pi 
A_3d^{p-1} + o(d^{p-1})\right)^{1/q}\right\}^{(p-1)/(p-3)}
\nonumber
\\
&& \times \left\{E_3^{-(p-3)/2}\alpha^{p-3}
\left(1 - \frac{p-3}{2}E_5\alpha^{p-3} + o(\alpha^{p-3})\right)\right\}^{2/(p-3)}
\nonumber
\\
&&\times\left\{
\pi^2 + 2\pi A_3d^{p-1} + o(d^{p-1})
\right\}^{{\{q(p-3)-(p-1)\}/(q(p-3))}}
\nonumber
\\
&=&
\left[(a_12^{(q+2)/q}A_1^{2/q} + a_2\pi^{2/q}) 
\right.
\\
&&
\left. \qquad \qquad + 
\left\{a_12^{(q+2)/q}A_1^{2/q}\left(A_4 + \frac{2}{q}\frac{A_2}{A_1}\right) 
+ \frac{2a_2\pi^{(2-q)/q}}{q}A_3\right\} d^{p-1} \right]
\nonumber
\\
&& \times E_3^{-1}\alpha^2\left\{1 - E_5\alpha^{p-3}+ o(\alpha^{p-3})\right\}
\nonumber
\\
&&\times \pi^{2{\{q(p-3)-(p-1)\}/(q(p-3))}}
\left\{1 + \frac{2\{q(p-3)-(p-1)\}}{(p-3)q\pi}A_3 d^{p-1} + o(d^{p-1})\right\}
\nonumber
\\
&=& \pi^{2\{q(p-3)-(p-1)\}/(q(p-3))}E_1E_3^{-1}\alpha^2 
\nonumber
\\
&& \qquad \qquad \times 
\left[
1 + \left\{\frac{E_2}{E_1} + E_4 - E_3^{(p-3)/2}E_5\right\}E_3^{-(p-3)/2} \alpha^{p-3} 
+ o(\alpha^{p-3})
\right].
\nonumber
\end{eqnarray}
This implies (1.5). Thus the proof is complete. \qed

\end{document}